\documentclass[reqno]{amsart}
\usepackage{hyperref}
\usepackage{latexsym}
\usepackage{graphicx}
\usepackage{amsmath}

\newtheorem{theorem}{Theorem}
\newtheorem{lemma}[theorem]{Lemma}
\newtheorem{corollary}[theorem]{Corollary}
\newtheorem{proposition}[theorem]{Proposition}
\newtheorem{example}[theorem]{Example}
\newcommand{\bt}{\begin{theorem}}
\newcommand{\et}{\end{theorem}}
\newcommand{\blem}{\begin{lemma}}
\newcommand{\elem}{\end{lemma}}
\newcommand{\bex}{\begin{example} \rm }
\newcommand{\eex}{\end{example}}
\newcommand{\bcor}{\begin{corollary}}
\newcommand{\ecor}{\end{corollary}}

\newcommand{\R}{\mathbb R}

\newcommand{\beq}{ \begin{equation} }
\newcommand{\eeq}{\end{equation} }
\newcommand{\bea}{\begin{eqnarray}}
\newcommand{\eea}{\end{eqnarray}}
\newcommand{\beas}{\begin{eqnarray*}}
\newcommand{\eeas}{\end{eqnarray*}}
\newcommand{\beqs}{\begin{equation*}}
\newcommand{\eeqs}{\end{equation*}}

\newcommand{\ben}{\begin{enumerate}}
\newcommand{\een}{\end{enumerate}}
\newcommand{\ba}{\begin{array}}
\newcommand{\ea}{\end{array}}

\newcommand{\eps}{\varepsilon}

\newcommand{\Fl}{\mathrm{Fl}}

\begin{document}

\title[Global Gronwall Estimates on Riemannian Manifolds]%
{Global Gronwall Estimates for Integral Curves on Riemannian Manifolds}

\author[M. Kunzinger]{Michael Kunzinger}
\address{Fakult\"at f\"ur Mathematik, Universit\"at  Wien\\
Nordbergstrasse 15\\ 1090 Wien\\ Austria}
\email{Michael.Kunzinger@univie.ac.at}
\urladdr{http://www.mat.univie.ac.at/\~{}mike/}
\author[H. Schichl]{Hermann Schichl}
\address{Fakult\"at f\"ur Mathematik, Universit\"at  Wien\\
Nordbergstrasse 15\\ 1090 Wien\\ Austria}
\email{Hermann.Schichl@univie.ac.at}
\urladdr{http://www.mat.univie.ac.at/\~{}herman/}
\author[R. Steinbauer]{Roland Steinbauer}
\address{Fakult\"at f\"ur Mathematik, Universit\"at Wien\\
Nordbergstrasse 15\\ 1090 Wien\\ Austria}
\email{Roland.Steinbauer@univie.ac.at}
\urladdr{http://www.mat.univie.ac.at/\~{}stein/}
\author[J. Vickers]{James A Vickers}
\address{School of Mathematics, University of Southampton\\ 
Highfield, Southampton SO {17 1}BJ, UK}
\email{J.A.Vickers@maths.soton.ac.uk}
\urladdr{http://www.maths.soton.ac.uk/staff/Vickers/}
\thanks{Partially supported by the Austrian Science Fund
Projects P16742-N04 and START-Project Y237-N13}

\begin{abstract}
We prove Gronwall-type estimates for the distance of
integral curves of smooth vector fields on a Riemannian manifold.
Such estimates are of central importance for all methods of solving ODEs in a verified way, 
i.e., with full control of roundoff errors. Our results may therefore
be seen as a prerequisite for the  generalization of such methods to the setting
of Riemannian manifolds.

\end{abstract}

\maketitle

\section{introduction}
Suppose that $X$ is a complete smooth vector field on $\R^n$, let $p_0$, $q_0 \in \R^n$
and denote by $p(t)$, $q(t)$ the integral curves of $X$ with initial values
$p_0$ resp.\ $q_0$. In the theory of ordinary differential equations it is a well 
known consequence of Gronwall's inequality  that in this situation we have
\begin{equation} \label{localgron}
|p(t)-q(t)| \le |p_0-q_0|e^{C_Tt} \qquad (t\in [0,T))
\end{equation}
with $C_T=\|DX\|_{L^\infty(K_T)}$ ($K_T$ some compact convex set containing the
integral curves $t\mapsto p(t)$ and $t\mapsto q(t)$) and $DX$ the Jacobian of $X$
(cf., e.g., \cite{dieudonne}, 10.5).

The aim of this paper is to derive estimates analogous to (\ref{localgron})
for integral curves of vector fields on Riemannian manifolds. Apart from
a purely analytical interest in this generalization, we note that 
Gronwall-type estimates play an essential role in the convergence analysis of 
numerical methods for solving ordinary differential equations (cf.\ \cite{stoer}).
Concerning notation and terminology from Riemannian geometry our basic references 
are \cite{ghl, klingenberg, oneill}. 

\section{Estimates}
The following proposition provides the main technical ingredient for the proofs of 
our Gronwall estimates. Here and in what follows, for $X\in {\mathfrak X}(M)$ (the space
of  smooth vector fields on $M$) we denote by $\nabla X$ its covariant differential and by 
$\|\nabla X(p)\|_g$ the mapping norm of $\nabla X(p): (T_pM,\|\,.\,\|_g) \to (T_pM,\|\,.\,\|_g),
\, Y_p \mapsto \nabla_{Y_p}X$. 
\begin{proposition} \label{main}
Let $[a,b]\ni \tau \mapsto c_0(\tau)=:c(0,\tau)$ be a smooth regular curve in a Riemannian manifold $(M,g)$, 
let $X\in {\mathfrak X}(M)$ and set $c(t,\tau) := \Fl^X_t c(0,\tau)$
where $Fl^X$ is the flow of $X$. 
Choose $T>0$ such that $\Fl^X$ is defined on $[0,T]\times c_0([a,b])$.
Then denoting by $l(t)$ the length of $\tau \mapsto c(t,\tau)$, we have
\begin{equation}
  \label{prop1equ}
  l(t) \leq l(0)e^{C_Tt} \qquad (t\in [0,T])
\end{equation}
where $C_T=\sup\{\|\nabla X(p)\|_g  : p\in c([0,T]\times [a,b])\}$.
\end{proposition}

\begin{proof}
Let $\tau \mapsto c(0,\tau)$ be parametrized by arclength, $\tau\in [0,l(0)]$.
Since $\Fl^X_t$ is a local diffeomorphism, $g(\partial_\tau c,\partial_\tau c)>0$ on $[0,T]\times [a,b]$.
Furthermore, since the Levi Civita connection $\nabla$ is torsion free, we have $\nabla_{\partial_t}c_\tau = 
\nabla_{\partial_\tau}c_t$, where $c_t = \partial_t c$, $c_\tau = \partial_\tau c$, see \cite{klingenberg},
1.8.14. Then
 \begin{eqnarray*}
  l(s)-l(0)
  &=&\int\limits_0^{s}\partial_t l(t)\,dt
  \,=\,\int\limits_0^{s}\partial_t\int\limits_0^{l(0)}\| c_\tau(t,\tau)\|_g\, d\tau\, dt\\
  &=&\int\limits_0^{s}\int\limits_0^{l(0)}
   \frac{\partial_t g(c_\tau(t,\tau),c_\tau(t,\tau))}{2 \| c_\tau(t,\tau)\|_g}\, d\tau\,dt
  \,=\,\int\limits_0^{s}\int\limits_0^{l(0)}
   \frac{g((\nabla_{\partial_t} c_\tau) (t,\tau),c_\tau(t,\tau))}{\| c_\tau(t,\tau)\|_g}\,d\tau\,dt \\
&=& 
\int\limits_0^{s}\int\limits_0^{l(0)}
   \frac{g((\nabla_{\partial_\tau} c_t) (t,\tau),c_\tau(t,\tau))}{\| c_\tau(t,\tau)\|_g}\,d\tau\,dt 
\le
\int\limits_0^{s}\int\limits_0^{l(0)} \|(\nabla_{\partial_\tau} c_t)(t,\tau)\|_g    \,d\tau\,dt\,,\\
&=&\int\limits_0^{s}\int\limits_0^{l(0)} \|\nabla_{c_\tau(t,\tau)} X\|_g    \,d\tau\,dt
\le C_T \int\limits_0^{s}\int\limits_0^{l(0)} \|c_\tau(t,\tau)\|_g    \,d\tau\,dt \\
&=& C_T \int_0^s l(t)\,dt\,.
\end{eqnarray*}
The claim now follows by applying Gronwall's inequality.
\end{proof}
We may utilize this proposition to prove our first main result:

\begin{theorem} \label{general}
Let $(M,g)$ be a connected smooth Riemannian manifold, $X\in {\mathfrak X}(M)$ a complete
vector field on $M$ and let $p_0,\, q_0 \in M$. Let $p(t) = \Fl^X_t(p_0)$, 
$q(t) = \Fl^X_t(q_0)$ and suppose that $C:=\sup_{p\in M}||\nabla X(p)||_g < \infty$.
Then
\begin{equation} \label{gronwallest}
d(p(t),q(t)) \le d(p_0,q_0) e^{Ct} \qquad (t\in [0,\infty))\,,
\end{equation}
where $d(p,q)$ denotes Riemannian distance.
\end{theorem}
\begin{proof}
For any given $\eps>0$, choose a piecewise smooth regular curve $\tau\mapsto c_0(\tau)=:c(0,\tau):
[0,1]\to M$ 
connecting $p_0$ and $q_0$ such that $d(p_0,q_0)>l(0)-\eps$. 
Using the notation of Proposition \ref{main} it follows that
$$
d(p(t),q(t)) \le l(t) \le l(0)e^{Ct} < (d(p_0,q_0)+\eps)e^{Ct}
$$
for $t\in [0,\infty)$. Since $\eps>0$ was arbitrary, the result follows.
\end{proof}

\bex \label{counterex}

\begin{figure}
\centerline{\resizebox{7cm}{!}{\includegraphics{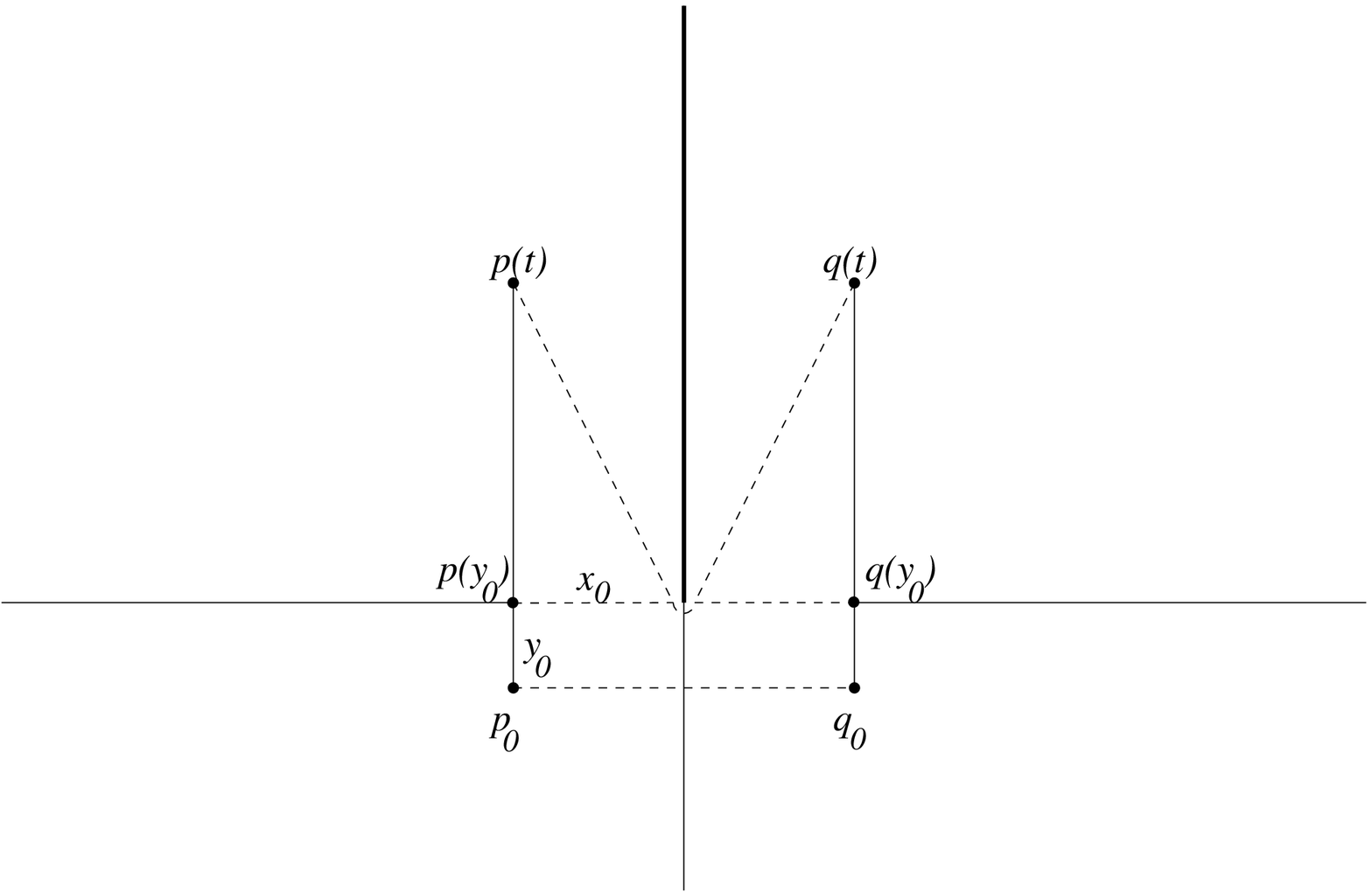}}}
\caption{}\label{noarmadillo}
\end{figure}

\noindent(i) In general, when neither $M$ nor $X$ is complete, the conclusion
of Theorem \ref{general} is no longer valid: 

Consider $M=\R^2\setminus \{(0,y)\mid y\ge 0\}$, endowed with the 
standard Euclidean metric. Let $X\equiv(0,1)$, $p_0=(-x_0,-y_0)$, and 
$q_0=(x_0,-y_0)$ ($x_0>0$, $y_0 \ge 0$) (cf.\ Figure \ref{noarmadillo}). 
Then $p(t)= (-x_0,-y_0+t)$, $q(t)=(x_0,-y_0+t)$ and 
$$
d(p(t),q(t))= \left\{ 
\begin{array}{cr}
2x_0 & \quad t\le y_0 \\
2\sqrt{x_0^2+(t-y_0)^2} &  \quad t>y_0
\end{array}
\right.
$$
On the other hand, $\nabla X=0$, so (\ref{gronwallest}) is violated for $t>y_0$,
i.e., as soon as the two trajectories are separated by the ``gap''
$\{(0,y)\mid y\ge 0\}$.

\noindent (ii) Replace $X$ in (i) by the complete vector field $(0,e^{-1/x^2+1})$
and set $x_0=1$, $y_0=0$. Then $C:=\|\nabla X\|_{L^\infty(\R^2)} = 3\sqrt{3/(2e)}$ and
$$
d(p(t),q(t)) = 2 \sqrt{1+t^2} \le d(p_0,q_0) e^{Ct} = 2 e^{Ct}
$$
for all $t\in [0,\infty)$, in accordance with Theorem \ref{general}.
\eex

The following result provides a sufficient condition for the validity of
a Gronwall estimate even if neither $M$ nor $X$ satisfies a completeness
assumption.

\bt \label{comp}
Let $(M,g)$ be a connected smooth Riemannian manifold, $X\in {\mathfrak X}(M)$ 
and let $p_0,\, q_0 \in M$. Let $p(t) = \Fl^X_t(p_0)$, 
$q(t) = \Fl^X_t(q_0)$ and suppose that there exists some relatively compact
submanifold $N$ of $M$ containing $p_0$, $q_0$ such that $d(p_0,q_0) = d_N(p_0,q_0)$. 
Fix $T>0$ such that $\Fl^X$ is defined on 
$[0,T]\times N$ and set $C_T:=\sup\{||\nabla X(p)||_g : p\in \Fl^X([0,T]\times N)\}$.
Then
\begin{equation} \label{gronwallfinite}
d(p(t),q(t)) \le d(p_0,q_0) e^{C_Tt} \qquad (t\in [0,T])\,.
\end{equation} 
\et
\begin{proof} 
As in the proof of Theorem \ref{general}, for any given $\eps>0$ we may choose a
piecewise smooth regular curve $\tau \mapsto c_0(\tau):[0,1]\to N$ from $p_0$ to $q_0$ such that
$d(p_0,q_0) = d_N(p_0,q_0) > l(0)-\eps$. The corresponding time evolutions
$c(t,\,.\,)$ of $c(0,\,.\,)=c_0$ then lie in $\Fl^X([0,T]\times N)$, so an application of
Proposition \ref{main} gives the result.
\end{proof}
\bex Clearly such a submanifold $N$ need not exist in general. As 
a simple example take $M=\R^2\setminus\{(0,0)\}$, $p_0=(-1,0)$, $q_0=(1,0)$.
In Example \ref{counterex}.(i) with $y_0>0$ the condition is obviously satisfied with
$N$ an open neighborhood of the straight line joining $p_0$, $q_0$ 
and the supremum of the maximal evolution times
of such $N$ under $\Fl^X$ is  $T=y_0$, coinciding with the maximal time-interval
of validity of (\ref{gronwallfinite}). On the other hand, if there is no
$N$ as in Theorem \ref{comp} then the conclusion in general
breaks down even for arbitrarily close initial points $p_0$, $q_0$:
if we set $y_0=0$ in Example \ref{counterex}.(i) 
then no matter how small $x_0$ (i.e., irrespective of the initial distance
of the trajectories) the estimate is not valid for any $T>0$. 
\eex
Finally, we single out some important special cases of Theorem \ref{comp}:
\bcor \label{cor}
Let $M$ be a connected geodesically complete Riemannian manifold,
$X\in {\mathfrak X}(M)$, and $p_0$, $q_0$, $p(t)$, $q(t)$ as above. 
Let $S$ be a minimizing geodesic segment connecting $p_0$, $q_0$
and choose some $T>0$ such that $\Fl^X$ is defined on 
$[0,T]\times S$. Then (\ref{gronwallfinite}) holds with 
$C_T=\sup\{||\nabla X(p)||_g\mid p\in \Fl^X([0,T]\times S)\}$.
In particular, if $X$ is complete then for any $T>0$ we have
\begin{equation} \label{gronwallest2}
d(p(t),q(t)) \le d(p_0,q_0) e^{C_Tt} \qquad (t\in [0,T])\,.
\end{equation}
\ecor
\begin{proof}
Choose for $N$ in Theorem \ref{comp} any relatively compact open neighborhood 
of $S$. The value of $C_T$ then follows by continuity.  
\end{proof}
In particular, for $M=\R^n$ with the standard Euclidean metric, Corollary 
\ref{cor} reproduces  (\ref{localgron}).

\end{document}